\documentclass[12pt,twoside,reqno]{amsart}
\usepackage[colorlinks=true,citecolor=blue]{hyperref}
\usepackage{mathptmx, amsmath, amssymb, amsfonts, amsthm, mathptmx, enumerate, color,mathrsfs}
\usepackage{graphicx}

\usepackage{multirow}
\usepackage{epstopdf}
\usepackage{multicol}
\usepackage{algorithm}
\usepackage{algorithmic}
\usepackage{epstopdf}

\newtheorem{theorem}{Theorem}[section]

\newtheorem{corollary}[theorem]{Corollary}

\theoremstyle{definition}
\newtheorem{definition}[theorem]{Definition}

\newtheorem{example}[theorem]{Example}

\numberwithin{equation}{section}

\def\Halmos{\mbox{\quad$\square$}}
\def\proof#1{\noindent {\it #1}}
\def\endproof{\vspace{0pt}}

\newcommand{\slim} {\mathop{\rm lim\,sup}}
\newcommand{\ilim} {\mathop{\rm lim\,inf}}
\def\R{\mathbb{R}}

\def\olR{{\overline{\R}}}

\def\S{\mathbb{S}}				

\def\ee{\varepsilon}

\def\h{\mathbf{I}}
\def\N{\mathbb{N}}			

\def\X{\mathbb{X}}

\def\A{\mathbb{A}}
\def\K{\mathbb{K}}
\def\U{\mathbb{U}}
\def\a{\alpha}				
\def\K{\mathbb{K}}			
\def\PS{\Pi}					
\def\PR{P}					

\begin{document}
\setcounter{page}{1}

\vspace*{2.0cm}
\title[Average-Cost MDPs with Infinite State and Action Sets]
{Average-Cost MDPs with Infinite State and Action Sets: New Sufficient Conditions for Optimality Inequalities and Equations}
\author[E.A.~Feinberg, P.O.~Kasyanov, L.S.~Paliichuk]{ Eugene~A.~Feinberg$^{1,*}$, Pavlo~O.~Kasyanov$^{2}$, Liliia~S.~Paliichuk$^2$}
\maketitle
\vspace*{-0.6cm}

\begin{center}
{\footnotesize

$^1$Department of Applied Mathematics and Statistics, Stony Brook University, Stony Brook, NY 11794-3600, USA\\
$^2$Institute for Applied System Analysis, National Technical University of Ukraine ``Igor Sikorsky Kyiv Polytechnic Institute'', Beresteiskyi Ave., 37, build 35, Kyiv, 03056 Ukraine

}\end{center}

\begin{center}
    {\it {This paper is dedicated to the memory of Uriel G. Rothblum}}
\end{center}

\vskip 4mm {\footnotesize \noindent {\bf Abstract.}
This paper studies discrete-time average-cost infinite-horizon Markov decision processes (MDPs) with Borel state and action sets. It introduces new sufficient conditions for { the} validity of optimality  inequalities and optimality equations for MDPs with weakly and setwise continuous transition probabilities. {These inequalities and equations imply the existence of deterministic optimal policies.}


  }

\renewcommand{\thefootnote}{}
\footnotetext{ $^*$Corresponding author.
\par
E-mail address: eugene.feinberg@stonybrook.edu (E.A.~Feinberg), kasyanov@i.ua (P.O.~Kasyanov), Lili262808@gmail.com (L.S.~Paliichuk).}


\section{Introduction}  This paper establishes sufficient conditions for the existence of deterministic optimal policies minimizing expected costs per unit time for infinite-horizon Markov Decision Processes with infinite state and action sets.  Such policies exist for problems with finite state and actions sets \cite{Bl,De,VS}, and deterministic policies were called stationary or randomized stationary in earlier publications.  However, if either the state space or the action space is infinite, optimal policies may not exist.  In particular, for countable-state MDPs with finite action sets, there are examples demonstrating nonexistence of optimal policies \cite{Ro, DY, Fe2}. For finite-state MDPs there are examples when optimal policies do not exist when action sets are compact, costs do not depend on actions, and one-step transition probabilities depend continuously on actions \cite{Ba1, Ch, DY}.

For finite-state MDPs with compact action sets, deterministic optimal policies exist in the following two cases: (i) all sets of transition probabilities have finite sets of extreme points \cite{Fe1}, (ii) the MDP is communicating \cite{Ba2}, that is, any state can be reached from any state.

For countable state MDPs, Sennott~\cite{Se1,Se2} proved the validity of optimality inequalities under conditions generalizing the communicating condition, and these inequalities imply the existence of deterministic optimal policies.   Cavazos-Cadena~\cite{Cad91} provides an example of a countable-state MDP, for which optimality inequalities hold, and the optimality equation does not.  Sch\"al~\cite{Sch93} extended these results to MDPs with possibly uncountable state sets and compact action sets by {considering} Assumption~B formulated below. { In view of Hern\'{a}ndez-Lerma and Lasserre~\cite[Theorem 5.4.6]{HLL96}, Assumptio~B is equivalent to the assumptions for the validity of the  optimality inequality in Sennott~\cite{Se1,Se2}.}  Feinberg et al.~\cite{FKZ12} and Feinberg and Kasyanov~\cite{MDPSetwise2021} extended Sch\"al's~\cite{Sch93} results to MDPs with possibly noncompact action sets and introduced a weaker Assumption~\underline{B}, which also implies an optimality inequality in a weaker form, which also implies the existence of deterministic optimal policies. \cite[Example 4.1]{MDPSetwise2021} demonstrates that Assumption~\underline{B} is indeed weaker than Assumption~B.   The results for noncompact action sets are important for inventory control~\cite{Ftut,FLe17}.

Sch\"al~\cite{Sch93} studied MPDs with weakly and setwise continuous transition probabilities.  Though weak continuity is more general than setwise continuity, MDPs with weakly continuous transition probabilities are not more general since continuity of costs and transition probabilities is assumed with respect to state-action pairs, while, for MDPs with setwise transition probabilities,  continuity of costs and transition probabilities is assumed only with respect to actions. Feinberg et al.~\cite{FKZ12} studied MDPs with weakly continuous transition probabilities, and Feinberg and Kasyanov~\cite{MDPSetwise2021} studied MDPs with setwise continuous transition probabilities.  Both models have important applications.  For example, MDPs with weakly continuous transition probabilities are used for partially observable MDPs \cite{FKZ14, POMDP}. An MDP with finite action sets and with arbitrary transition probabilities and arbitrary costs is an example of an MDP with setwise continuous transition probabilities. Hern\'andez-Lerma~\cite{HL91} studied MDPs with setwise continuous transition probabilities with possibly noncompact action sets under Assumption~B, but the optimality equation for discounted MDPs, which was used in {the} proofs, was formulated there without a proof, and the only proof known to the authors follows from the optimal selection theorem proved later in \cite{MDPSetwise2021}.

In this paper we introduce sufficient conditions, which are weaker than Assumption B, and which lead to the same {or weaker} conclusions on the validity of optimality inequalities and optimality equations as Assumption~B.  There is a significant literature on MDPs with average costs per unit time, which includes three surveys~\cite{Ari,Bor,Se2}.  Recently Guo et al.~\cite{GHZ} established new conditions for the existence of optimal policies for non-stationary MDPs.

\section{Preliminaries}

{Let $\olR:=\mathbb{R}\cup\{+\infty\},$ $\mathbb{N}^*:=\{0,1,\ldots\}=\mathbb{N}\cup\{0\}.$} Consider a discrete-time MDP with a state space $\X,$ an action space $\A,$ one-step
costs $c,$ and transition probabilities $q.$ Assume that $\X$ and $\A$ are Borel subsets
of Polish (complete separable metric) spaces. 
Let  $c(x,a):\X\times\A\mapsto\olR$ be the one-step cost and $q(B|x,a)$ be the transition kernel
representing the probability that the next state is in $B\in\mathcal{B}(\X),$
given that the action $a$ is chosen at the state $x.$  The cost function $c$ is assumed to be measurable and bounded below.

The decision process proceeds as follows: at each time epoch
$t=0,1,\dots,$ the current state of the system, $x,$ is observed.
A decision-maker chooses an action $a,$ the cost $c(x,a)$ is
accrued, and the system moves to the next state according to
$q(\,\cdot\,|x,a).$ Let $H_t = (\X\times\A)^{t}\times\X$
be the set of histories for $t=0,1,\dots\ .$ A
(randomized) decision rule at period $t=0,1,\dots$ is a regular transition probability
$\pi_t : H_t\mapsto \A,$ 
that is, (i) $\pi_t(\,\cdot\,|h_t)$ is a probability distribution on $\A,$
where $h_t=(x_0,a_0,x_1,\dots,a_{t-1},x_t),$ and (ii) for any measurable subset
$B\subset \A,$ the function $\pi_t(B|\cdot)$ is measurable on $H_t.$
A policy $\pi$ is a sequence $(\pi_0,\pi_1,\dots)$ of decision rules.  Let $\PS$ be the set of all policies.
A policy $\pi$ is called non-randomized if each probability measure $\pi_t(\,\cdot\,|h_t)$ is
concentrated at one point. A non-randomized policy is called deterministic if all decisions depend
only on the current state.

The Ionescu Tulcea theorem implies that an initial state $x$ and a policy $\pi$ define a unique
probability $\PR_{x}^{\pi}$ on the set of all trajectories $\mathbb{H}_{\infty}=(\X\times\A)^{\infty}$ endowed with the product of $\sigma$-fields defined by Borel $\sigma$-fields of $\X$ and $\A;$ see Bertsekas and Shreve~\cite[pp. 140--141]{BS96} or Hern\'{a}ndez-Lerma and Lasserre~\cite[p. 178]{HLL96}.
Let $\mathbb{E}_{x}^{\pi}$ be an expectation w.r.t. $\PR_{x}^{\pi}.$

For a finite-horizon $N\in\N^*:=\{1,2,\ldots\},$ let us define the expected total discounted costs,
\begin{align}\label{eqn:sec_model def:finite total disc cost}
    v_{N,\a}^{\pi} (x):= \mathbb{E}_{x}^{\pi} \sum_{t=0}^{N-1}
    \alpha^{t} c(x_t,a_t), \;\; x\in\X,
\end{align}
where $\alpha\in [0,1]$ is the discount factor.
When $N=\infty$ and $\alpha\in [0,1),$
equation (\ref{eqn:sec_model def:finite total disc cost}) defines an
infinite-horizon expected total discounted cost denoted by $v_{\a}^{\pi}(x).$  {We always assume that $\alpha\in [0,1)$ when $N=\infty.$  We observe that the expectation in \eqref{eqn:sec_model def:finite total disc cost} is well-defined in the following two cases: (i) $N<\infty,$  (ii) $N=\infty$ and $\alpha\in [0,1).$ This is true because the sum in (\ref{eqn:sec_model def:finite total disc cost}) is a bounded below measurable function since the function $c$ is  bounded below and measurable.}

Let $v_\alpha (x):=\inf_{\pi\in\PS} v_\a^\pi(x),$ $x\in\X.$  A policy $\pi$ is called optimal for the discount factor $\a$ if $v^\pi_\a(x)=v_\a(x)$ for all $x\in\X.$

The \emph{average cost per unit time} is defined as
\begin{align*}
    w^{\pi}(x):=\limsup_{N\to \infty} \frac{1}{N}v_{N,1}^{\pi} (x), \;\; x\in\X.
\end{align*}
Define the optimal value function $w(x):=\inf_{\pi\in\Pi} w^{\pi} (x),$  $x\in\X.$
A policy $\pi$ is called average-cost optimal if $w^{\pi}(x)=w(x)$ for all $x\in\X.$

We remark that in general action sets may depend on current states, and usually the state-dependent sets $A(x)$ are considered for all $x\in\X$.  In our problem formulations $A(x)=\A$ for all $x\in\X.$  This problem formulation is simpler than a  formulation with the sets $A(x),$ and these two problem formulations are equivalent because we allow that $c(x,a)=+\infty$ for some $(x,a)\in\X\times\A$  and can set $A(x)=\{a\in\A: c(x,a)<+\infty\}.$  For  a  formulation with the sets $A(x),$ one may define $c(x,a)=+\infty$ when $a\in\A \setminus A(x)$ and use the action {set} $\A$ instead of $A(x).$

To establish the existence of  average-cost optimal policies  for problems with compact action sets,
Sch\"{a}l~\cite{Sch93} considered two continuity Assumptions~W and S for problems with weakly and setwise continuous transition probabilities, respectively.
For setwise continuous transition probabilities,
Hern\'{a}ndez-Lerma~\cite{HL91} generalized
Assumption~S to Assumption~S*  to cover MDPs with possibly noncompact action sets.
For the similar purpose, when transition probabilities are weakly continuous,
Feinberg~et~al.~\cite{FKZ12} generalized Assumption~\textbf{W} to Assumption~W*.

We recall that a function $f:\U\mapsto \olR$ defined on a metric space $\U$  is
called inf-compact (on $\U$), if for every $\lambda\in\R$ the level set $\{u\in\U :f(u)\leq \lambda \}$
is compact. A { measurable} subset of a metric space is also a metric space with respect to the same metric.  For $U\subset \U,$ if the domain of $f$ is narrowed to $U,$ then this function is  called the restriction of $f$ to $U.$
\begin{definition}[{Feinberg et al.~\cite[Definition 1.1]{FKZ13}, Feinberg~\cite[Definition 2.1]{Ftut}}]
\label{def:k inf compact}
	A function $f:\X\times\A\mapsto \olR$ is called $\K$-inf-compact, if for
	every nonempty compact subset $\mathcal{K}$ of $\X$ the restriction of $f$ to $\mathcal{K}\times\A$ is an inf-compact function.
\end{definition}

\noindent
\textbf{Assumption W*} ({\rm Feinberg~et~al.~\cite{FKZ12, POMDP}}, Feinberg and Lewis~\cite{FLe17}, or Feinberg~\cite{Ftut})\textbf{.}

(i) the function $c$ is $\K$-inf-compact; 

(ii) the transition probability $q(\,\cdot\,|x,a)$ is weakly continuous in
  $(x,a)\in \X\times\A.$

\vspace{0.1in}

\noindent
\textbf{Assumption S*} ({\rm Hern\'{a}ndez-Lerma~\cite[Assumption 2.1]{HL91}} or Feinberg and Kasyanov~\cite{MDPSetwise2021})

(i) the function $c(x,a)$ is inf-compact in $a\in \A$ for each $x\in\X;$

(ii) the transition probability $q(\,\cdot\,|x,a)$ is setwise continuous in
  $a\in \A$ for each $x\in\X .$
\vspace{0.1in}


Let
\begin{align}\label{defmauaw}
	\begin{split}
  		m_{\alpha}: = \underset{x\in\X}{\inf} v_{\alpha}(x), & \quad
  		u_{\alpha}(x): = v_{\alpha}(x) - m_{\alpha}, \\
  		\underline{w}: = \underset{\alpha\uparrow1}{\liminf}(1-\alpha)m_{\alpha}, & \quad
  		\bar{w}: = \underset{\alpha\uparrow1}{\limsup}(1-\alpha)m_{\alpha} .
	\end{split}
\end{align}
The function $u_\alpha$ is called the discounted relative value function.  If either Assumption~W* or Assumption~S* holds, let us consider the following assumption.

\vspace{0.1in}
\noindent
\textbf{Assumption B.} (Sch\"al~\cite{Sch93}).
(i) $w^{*} := \inf_{x\in\X} w (x)< +\infty;$ and
(ii) $\underset{\alpha\in [0,1)}{\sup} u_{\alpha}(x) < +\infty,$ $x\in\X.$

\vspace{0.1in}

{ We recall that $\underset{\alpha\in [0,1)}{\sup} u_{\alpha}(x) < +\infty$ if and only if $\slim\limits_{\alpha \uparrow
1}u_{\alpha}(x)<+\infty;$}  \cite[Lemma 5]{FKZ12}.   As follows from Sch\"al~\cite[Lemma 1.2(a)]{Sch93},
Assumption~\textbf{B}(i) implies that $m_\alpha< +\infty$
for all $\alpha\in [0,1).$  Thus, all the quantities in
\eqref{defmauaw} are defined.  

It is known \cite[Theorem 1]{FKZ12} that, if a deterministic policy $\phi$ satisfies the weakened average-cost optimality inequality (WACOI):
\begin{align}c(x, \phi(x)) + \int_{\X} u(y) q(dy \mid x, \phi(x))   &\leq \overline{w} + u(x),
\qquad x \in \X, \label{eqn:ACOIOver} \end{align}
for some nonnegative measurable function $u: \X \to {\mathbb R},$ then the deterministic policy $\phi$ is average-cost optimal, and
\begin{equation}\label{eq:7121} w(x) = w^{\phi}(x) = \lim_{\alpha \uparrow 1} (1 - \alpha) v_{\alpha}(x) = \overline{w} = w^*, \quad x \in \X. \end{equation}

{
Let us consider the following assumption.}

\vspace{0.1in}
\noindent
\textbf{Assumption~\underline{B}.} {(Feinberg et al. \cite{FKZ12})}.
(i) Assumption~B(i)
 holds, and (ii) \ $\ilim\limits_{\alpha \uparrow
1}u_{\alpha}(x)<+\infty$ for all $x\in \X$.
\vspace{0.1in}

Assumption~\underline{B}(ii) is weaker than
 Assumption~B(ii); see \cite[Example~4.1]{MDPSetwise2021}.
If a deterministic policy $\phi$ satisfies the average-cost optimality inequality (ACOI):
\begin{align} c(x, \phi(x)) + \int_{\X} u(y) q(dy \mid x, \phi(x)) &\leq  \underline{w} + u(x),
\qquad x \in \X, \label{eqn:ACOI} \end{align}
for some nonnegative measurable function $u: \X \to {\mathbb R},$ which is a stronger version of \eqref{eqn:ACOIOver} because $\underline{w} \leq \overline{w}$ always holds, then, according to \cite{Sch93}, the deterministic policy $\phi$ is average-cost optimal, and in addition to \eqref{eq:7121}, it follows that $\underline{w} = \overline{w}.$ A nonnegative measurable function $u(x)$ satisfying inequality \eqref{eqn:ACOI} with some deterministic policy $\phi$ is
called an average-cost relative value function.
``Boundedness'' Assumption~\underline{B} on the
function $u_\alpha$, which is weaker than boundedness
Assumption~ B, and either Assumption W* or Assumption S* lead
to the validity of WACOI \eqref{eqn:ACOIOver} and
the existence of optimal deterministic policies \cite[Theorem 3]{FKZ12} and \cite[Theorem 3.3]{MDPSetwise2021}. Stronger results, namely, the validity of ACOI \eqref{eqn:ACOI} hold if
Assumption~B holds instead of Assumption~\underline{B}; see \cite[Theorem 4]{FKZ12} and \cite{HL91}.

We recall that $\alpha\in [0,1)$ {for infinite-horizon problems,} and everywhere in this paper, if we consider a discount factor   $\alpha_n,$ we assume that $\alpha_n\in [0,1).$

\section{Main Results}

We study MDPs either with weakly continuous transition probabilities satisfying Assumption~W* or with setwise continuous transition probabilities satisfying Assumption~S*. In either case, {each of the Assumptions~B or \underline{B}} imply the validity of optimality inequalities and the existence of deterministic optimal policies~\cite{MDPSetwise2021,FKZ12}.  However, the results are stronger under Assumption~B.  In addition, under additional conditions, Assumption~B implies the validity of the optimality equation~\cite{FKL18a}.
We prove in Corollaries~\ref{corACOIwm} and \ref{ceqACEstr} that the results on the validity of optimality inequalities and optimality equations, that were established under Assumption~B, hold under more general assumptions introduced in this section.

Theorems~\ref{thm:ACOI:w} and \ref{thm:ACOI:s} state the validity of
 WACOI~\eqref{eqn:ACOIOver} under Assumption~W*
or Assumption~S* and under essentially weakened version of Assumption~\textbf{B}. For this purpose for an arbitrary fixed sequence $\alpha_n \uparrow 1$ we set:
\begin{equation}\label{eq:omega}
    \underline{w}_{\{\alpha_n\}:}=\ilim\limits_{n\to\infty}(1-\alpha_n)m_{\alpha_n},\quad\overline{w}_{\{\alpha_n\}}:=\slim\limits_{n\to\infty}(1-\alpha_n)m_{\alpha_n}.
\end{equation}
According to
Sch\"al \cite[Lemma 1.2]{Sch93}, Assumption~B(i) implies
\begin{equation}\label{eq:schal} 0\le \underline{w}\le  \underline{w}_{\{\alpha_n\}}\le\overline{w}_{\{\alpha_n\}}\le \overline{w}\le w^*<
+\infty.
\end{equation}
The following theorem formulates  average-cost optimality inequality \eqref{eq7111} in a different form than ACOI~\eqref{eqn:ACOI} and WACOI~\eqref{eqn:ACOIOver}.

\begin{theorem}\label{Prop1}
Let Assumption~B(i) hold and $\{ \a_n\uparrow 1 \}_{n\in\N^*}$ be an arbitrary fixed sequence. If there exists a measurable
function $u:\X\to [0,+\infty)$ and a deterministic policy $\phi$ such
that
\begin{equation}\label{eq7111}
c(x, \phi(x))+\int_\X  u(y)q(dy|x,
\phi(x))\le \overline{w}_{\{\alpha_n\}}+u(x),\quad x\in \X,
\end{equation}
then $\phi$ is average-cost optimal,
\begin{equation}\label{eq:7121a}
w(x)=w^{\phi}(x)=\slim\limits_{n\to\infty}(1-\alpha_n)v_{\alpha_n}(x)=\slim\limits_{\alpha\uparrow
1}(1-\alpha)v_{\alpha}(x)=\overline{w}=\overline{w}_{\{\alpha_n\}}=w^*,
\end{equation}
for each $ x\in \X, $ and WACOI \eqref{eqn:ACOIOver} hold for the same policy $\phi$ and function $u$ as in \eqref{eq7111}.
\end{theorem}

We remark that for $ \alpha_n \uparrow 1,$ if  $(1-\alpha_n)m_{\alpha_n}\to \overline{w},$ then \eqref{eq7111} coincides with WACOI \eqref{eqn:ACOIOver}, which is already stated in Theorem~\ref{Prop1}, and, if $(1-\alpha_n)m_{\alpha_n}\to \underline{w},$ then \eqref{eq7111} coincides with ACOI \eqref{eqn:ACOI}, which is an additional property; see Corollary~\ref{corACOIwm}.

\proof{Proof of Theorem~\ref{Prop1}.} Similarly to  Feinberg et al. \cite[Theorem~1]{FKZ12}, since $u$ is nonnegative,
by iterating (\ref{eq7111}) we obtain
\[
v_{n,1}^\phi(x)\le n\overline{w}_{\{\alpha_n\}}+u(x), \quad n\ge 1,\ x\in \X.
\]
Therefore, after dividing the last inequality by $n$ and setting
$n\to\infty$, we have
\begin{equation}\label{eq***}
w^*\le w(x)\le w^{\phi}(x)\le
\overline{w}_{\{\alpha_n\}},\quad  x\in \X,
\end{equation}
where the first and the second inequalities follow from the
definitions of $w$ and $w^*$ respectively. Since $w^*\le \overline{w}_{\{\alpha_n\}},$ inequalities (\ref{eq:schal}) imply that for all $\pi\in\Pi$ and for all $x\in\X$
\[
w^*=\overline{w}=\overline{w}_{\{\alpha_n\}} \le \slim\limits_{n\to\infty}(1-\alpha_n)v_{\alpha_n}(x)\le \slim\limits_{\alpha\uparrow
1}(1-\alpha)v_{\alpha}(x)\le\slim\limits_{\alpha\uparrow
1}(1-\alpha)v_{\alpha}^{\pi}(x)\le w^{\pi}(x),
\]
where the last inequality follows from the Tauberian theorem. Finally, we obtain that
\begin{equation}\label{eq3.5}
\begin{aligned}
w^*=&\overline{w}=\overline{w}_{\{\alpha_n\}} \le \slim\limits_{n\to\infty}(1-\alpha_n)v_{\alpha_n}(x)\le \slim\limits_{\alpha\uparrow
1}(1-\alpha)v_{\alpha}(x)\\ &\le \inf\limits_{\pi\in
\Pi}w^{\pi}(x)=w(x)\le w^{\phi}(x) \le\overline{w}_{\{\alpha_n\}},
\end{aligned}
\end{equation}
for each $x\in\X,$ where the last inequality 
follows from (\ref{eq***}).  Thus, all the inequalities in
(\ref{eq3.5}) are equalities. \Halmos\endproof

For a sequence $\{ \a_n\uparrow 1 \}_{n\in\N^*} $ of discount factors, consider the following assumption.

\textbf{Assumption ${\rm   \underline{B}_{\{\alpha_n\}}}$.} (i)
Assumption~B(i)
 holds, and (ii) for a sequence $\{ \a_n\uparrow 1 \}_{n\in\N^*} $ of discount factors, the inequality $\ilim\limits_{n\to \infty}u_{\alpha_n}(x)<{+}\infty$ holds for all $x\in \X$.

 Assumption~B is equivalent to the statement that Assumption~${\rm  \underline{B}_{\{\alpha_n\}}}$ holds for an arbitrary sequence $\{\alpha_n \uparrow 1\}_{n \in \N^*}$  because Assumption B obviously implies this statement and, conversely, by contradiction, if Assumption B does not hold, then $\limsup_{n\to\infty}u_{\alpha_n}(x)\to +\infty$ for some sequence $\{\alpha_n \uparrow 1\}_{n \in \N^*},$ and  Assumption~${\rm  \underline{B}_{\{\beta_n\}}}$ does not hold for the subsequence $\{\beta_n\}_{n \in \N^*}$ of the sequence $\{\alpha_n\}_{n \in \N^*}$ such that $\lim_{n\to\infty}u_{\beta_n}(x)\to +\infty$. The existence of a sequence  $\{\alpha_n \uparrow 1\}_{n \in \N^*}$ satisfying Assumption~${\rm  \underline{B}_{\{\alpha_n\}}}$ implies Assumption~\underline{B}.
Moreover, we note that $\liminf_{n\to \infty, y\to x} u_{\alpha_n}(y)$ is the least upper bound of the set of all $\lambda \in \mathbb{R}_+$ such that there exist $m \in \mathbb{N}$ and a neighborhood $V(x)$ of $x$ satisfying
\[
\lambda \le \inf\{u_{\alpha_n}(y) : n \geq m,\, y \in V(x) \}.
\]
This holds because
\[
\liminf_{n\to \infty, y\to x} u_{\alpha_n}(y) = \sup_{V(x),\, m}\ \  \inf_{y \in V(x),\, n \ge m} u_{\alpha_n}(y),
\]
where the supremum is taken over all neighborhoods $V(x)$ of $x$ and $m =1,2,\ldots.$

\begin{theorem}\label{thm:ACOI:w}
Let Assumptions~W* holds and let Assumption ${\rm  \underline{B}_{\{\alpha_n\}}}$ hold for a sequence $\{ \a_n\uparrow 1 \}_{n\in\N^*}.$  Let
\begin{align}
  	u(x) := \liminf_{n\to \infty,y\rightarrow x} u_{\a_n}(y) , \qquad x\in\X .
  	\label{EQN1}
\end{align}
Then there exists a  deterministic  policy $\phi$ satisfying WACOI \eqref{eqn:ACOIOver} with the function $u$ defined in
\eqref{EQN1}.  Therefore, $\phi$ is a deterministic average-cost optimal policy.  In addition, the function $u$ is lower semi-continuous, and equalities \eqref{eq:7121a}
hold.
\end{theorem}

According {to} definition \eqref{EQN1} the function $u$ depends on the sequence $\{ \a_n\uparrow 1 \}_{n\in\N^*}. $ We do not write this dependence explicitly. A natural question, which we do not study in this paper, is under which conditions functions $u$ defined in \eqref{EQN1} coincide for two sequences of discount factors converging to 1.

{Note that the following properties take place in Example~4.1  from \cite{MDPSetwise2021}:} (a) Assumption~B does not hold; (b)  Assumptions~W* and ${\rm  \underline{B}_{\{\alpha_n\}}}$ for some sequence $\{ \a_n\uparrow 1 \}_{n\in\N^*} $ hold; (c) $\underline{w}=\overline{w}$ and, therefore,  there exists a  deterministic  policy $\phi$ satisfying ACOI \eqref{eq7111}, \eqref{eqn:ACOI} with the function $u$ defined in \eqref{EQN1}.

Let Assumption~${\rm  \underline{B}_{\{\alpha_n\}}}$ hold for a sequence $\{ \a_n\uparrow 1 \}_{n\in\N^*}.$ We define the following nonnegative functions on $\X$:
\begin{equation}\label{eq:defuuU}
U_m(x)  =  \inf\limits_{\small n\ge m}u_{\alpha_n}(x), \quad \underline{u}_m(x)  =
\ilim\limits_{y\to x}U_m(y),\ \quad\quad m=1,2,\ldots,\,\,
x\in \X.
\end{equation}
 Observe that all the three defined functions take finite
values at $x\in \X.$ Indeed,
\begin{equation}\label{eqnew9}
\underline{u}_m(x)\le U_m(x)\le \sup_{m=1,2,\ldots
}\inf\limits_{\small n\ge m}
u_{\alpha_n}(x)=\ilim\limits_{n\to\infty}u_{\alpha_n}(x)<{+}\infty,\quad m=1,2,\ldots, \ x\in\X,
\end{equation}
where the first two inequalities follow from the definitions of
$\underline{u}_m$ and $U_m$ respectively, and the last
inequality follows from Assumption ${\rm \underline{B}_{\{\alpha_n\}}}$. For
$x\in \X$
\begin{equation}\begin{split}\label{eq:ubounded}
u({x})=\sup\limits_{\small m=1,2,\ldots,\ R>0}\left[
\inf\limits_{\small n\ge m, \ y\in
B_{R}(x)}u_{\alpha_n}(y)\right]=\sup\limits_{m=1,2,\ldots}\
\sup\limits_{R>0}\ \inf\limits_{y\in B_R(x)} \
\inf\limits_{n\ge m}u_{\alpha_n}(y)\\ =
\sup\limits_{m=1,2,\ldots}\ \sup\limits_{R>0}\ \inf\limits_{y\in
B_R(x)} \ U_{m}(y) =\sup\limits_{m=1,2,\ldots}\
\ilim\limits_{y\to x}\ U_{m}(y)=\sup\limits_{m=1,2,\ldots}\underline{u}_{m}(x)<{+}\infty, \end{split}
\end{equation}
where $B_R(x)=\{y\in \X \,: \,    \rho(y,x)<R\}$, the first equality
is (\ref{EQN1}), the second equality follows from the properties of
infima, the third and the fifth equalities follow from
(\ref{eq:defuuU}), the fourth equality follows from the definition
of $\limsup$,  and the last inequality follows from (\ref{eqnew9}).
In view of (\ref{eq:defuuU}), the functions $U_m(x)$ and
$\underline{u}_m(x)$ are nondecreasing in $m$.  Therefore,
in view of (\ref{eq:ubounded}),
\begin{equation}\label{eq5.5}
u(x)=\lim\limits_{m\to\infty
}{\underline{u}_{m}}(x),\qquad\qquad x\in\X.
\end{equation}
 We
also set for $u$ from (\ref{eq5.5})
\begin{equation}\label{defsetA*}
 A^*(x):=\left\{a\in A(x)\, : \,c(x,a)+\int_\X  u(y)q(dy|x,a) \le
\overline{w}_{\{\alpha_n\}}+u(x)  \right\}, \
x\in\X.
\end{equation}

\proof{Proof of Theorem~\ref{thm:ACOI:w}.}
By replacing $\alpha\in[0,1)$ with $\alpha\in\{\alpha_n\}_n$ in Lemma~6 from Feinberg et al. \cite{FKZ12}, we obtain that the functions $u,u_{\alpha_m},$ and $
\underline{u}_m:\X\to \mathbb{R}_+,$ $m=1,2,\ldots,$ are
lower semi-continuous on $\X.$

{Let us prove} that $u$ satisfies \eqref{eq7111}. For this purpose, let us fix an arbitrary $\varepsilon^*>0$. Since $
\overline{w}_{\{\alpha_n\}}=\slim\limits_{n\to \infty}(1-\alpha_n)m_{\alpha_n}, $
there exists $n_0\in [0,1)$ such that
\begin{equation}\label{eq:l32}
\overline{w}_{\{\alpha_n\}} +\varepsilon^* > (1-\alpha_n)m_{\alpha_n},
\quad n=n_0,n_0+1,\ldots.
\end{equation}
Our next goal is to prove the inequality
\begin{equation}\label{eq:l33}
\overline{w}_{\{\alpha_n\}} +\varepsilon^*+u(x)\ge\min\limits_{a\in A(x)}
\left[c(x,a)+\alpha_m\int_\X  {\underline{u}_{m}}(y)q(dy|
x,a)\right],\quad x\in\X,\ m\ge n_0.
\end{equation}
Indeed, by
\begin{equation}\label{eq12}
(1-\alpha)m_{\alpha}+u_{\alpha}(x)=\min\limits_{a\in
A(x)}\left[c(x,a)+\alpha \int_\X  u_{\alpha}(y)q(dy|x,a)\right],
\quad x\in \X.
\end{equation}
and (\ref{eq:l32}) for every
$n,m\ge n_0$, such that $n\ge m$, and for
every $x\in \X$
\[
\overline{w}_{\{\alpha_n\}}+\varepsilon^*+u_{\alpha_n}(x)>
(1-\alpha_n)m_{\alpha_n}+u_{\alpha_n}(x)=\min\limits_{a\in
A(x)}\left[c(x,a)+\alpha_n\int_\X  u_{\alpha_n}(y)q(dy|x,a)\right]
\]
\[
\ge \min\limits_{a\in A(x)}\left[c(x,a)+\alpha_m\int_\X
U_{m}(y)q(dy|x,a)\right]
\]
because $\alpha_n\ge \alpha_m$ since $\alpha_n\uparrow 1.$
As the right-hand side does not depend on $n\ge m,$ we have
for all $x\in\X$ and for all $\alpha\in [\alpha_0,1)$
\[
\overline{w}_{\{\alpha_n\}}+\varepsilon^*+U_m(x)=
\inf\limits_{n\ge m }\left[\overline{w}_{\{\alpha_n\}}+\varepsilon^*+u_{\alpha_n}(x)
\right]
 \ge\min\limits_{a\in A(x)}\left[c(x,a)+\alpha_m\int_\X
U_{m}(y)q(dy|x,a)\right]
\]
\[
\ge\min\limits_{a\in A(x)}\left[c(x,a)+\alpha_m\int_\X
\underline{u}_{m}(y)q(dy|x,a)\right]= \min\limits_{a\in
A(x)}\eta_{{\underline{u}_{m}}}^{\alpha_m}(x,a),
\]
where
\[
\eta_{{\underline{u}_{m}}}^{\alpha_m}(x,a):= c(x,a)+\alpha_m\int_\X
\underline{u}_{m}(y)q(dy|x,a).
\]
By { Feinberg et al. \cite[Lemma~3]{FKZ12},} the function $x\mapsto \min\limits_{a\in
A(x)}\eta_{{\underline{u}_{m}}}^{\alpha_m}(x,a)$ is lower
semi-continuous on $\X$. Thus,
\[
\ilim\limits_{y\to x}\min\limits_{a\in
A(y)}\eta_{{\underline{u}_{m}}}^{\alpha_m}(y,a)\ge
\min\limits_{a\in A(x)}\eta_{{\underline{u}_{m}}}^{\alpha_m}(x,a),
\qquad x\in\X,\ m=1,2,\ldots.
\]
and, as, by definition (\ref{eq:defuuU}),
$\underline{u}_{m}(x)=\ilim\limits_{y\to x}{U}_{m}(y)$, we
finally obtain
\begin{equation}\label{eqalmweu}
\overline{w}_{\{\alpha_n\}}+\varepsilon^*+\underline{u}_m(x)\ge\min\limits_{a\in
A(x)}\eta_{\underline{u}_{m}}^{\alpha_m}(x,a),\quad
 x\in\X, m\ge n_0.
\end{equation}
Since $u(x)=\sup\limits_{m\ge n_0}\underline{u}_m(x)$ for all $x\in\X$,
(\ref{eqalmweu}) yields (\ref{eq:l33}).

To complete the proof of the theorem, we fix an arbitrary $x\in\X$.
By { Feinberg et al. \cite[Lemma~3]{FKZ12},} for any $m=1,2,\ldots$ there exists
$a_m\in A(x)$ such that $ \min\limits_{a\in
A(x)}\eta_{{\underline{u}_{m}}}^{\alpha_m}(x,a)=
\eta_{{\underline{u}_{m}}}^{\alpha_m}(x,a_{m}).$
 Since $\underline{u}_m\ge 0$, for
$m \ge n_0$ the inequality  (\ref{eq:l33}) can be
continued as
\begin{equation}\label{eq:5.14}
\overline{w}_{\{\alpha_n\}}+\varepsilon^*+u(x)\ge
\eta_{{\underline{u}_{m}}}^{\alpha_m}(x,a_m)\ge
c(x,a_m).
\end{equation}
Thus, for all $m\ge n_0$
\[
a_m\in
\mathcal{D}_{\eta_{\underline{u}_{m}}^{m}(x,\cdot)}(\overline{w}_{\{\alpha_n\}}+\varepsilon^*+u(x))\subset
\mathcal{D}_{c(x,\cdot)}(\overline{w}_{\{\alpha_n\}}+\varepsilon^*+u(x))\subset
A(x),
\]
where $\mathcal{D}_f(\lambda)=\{y\in U \, : \,  f(y)\le
\lambda\}$ is the level set.
By { Feinberg et al. \cite[Lemma~2]{FKZ12},} the set
$\mathcal{D}_{c(x,\cdot)}(\overline{w}_{\{\alpha_n\}}+\varepsilon^*+u(x))$ is
compact. Thus, there is a  subsequence $\{\alpha_m\}_{m} \subset \{\alpha_n\}_{n\ge 1}$ such that the sequence $\{a_{\alpha_m}\}_{m}$ converges and
$a_*:=\lim_{m} a_{\alpha_m}\in A(x)$.

Consider a subsequence $\{\alpha_m\}_m\subset \{\alpha_n\}_{n\ge1}$ such that $a_{\alpha_m}\to
a_*$ for some $a_*\in A(x).$ Due to Fatou's lemma for weakly converging probabilities \cite{FKL18},
\begin{equation}\label{eq:5.15} \ilim\limits_{m\to +\infty}
\alpha_{m}\int_\X
{\underline{u}_{m}}(y)q(dy|x,a_{m})\ge
\int_\X u(y)q(dy|x,a_*).
\end{equation}

Since the function $c$ is lower semi-continuous, (\ref{eq:5.14}) and
(\ref{eq:5.15}) imply
\[\overline{w}_{\{\alpha_n\}}+\varepsilon^*+u(x)\ge \limsup\limits_{n\to\infty}
\eta_{\underline{u}_{\alpha_n}}^{\alpha_n}(x,a_{\alpha_n})\ge
c(x,a_*)+\int_\X u(y)q(dy|x,a_*)\ge\min_{a\in A(x)} \eta_u^1(x,a).
\]
Since $\overline{w}_{\{\alpha_n\}}+\varepsilon^*+u(x)\ge\min_{a\in A(x)}
\eta_u^1(x,a)$ for all $\varepsilon^*>0$, this is also true when
$\varepsilon^*=0$.

 The Arsenin-Kunugui theorem
implies the existence of a deterministic policy $\phi$ such that
$\phi(x)\in A^*(x)$ for all $x\in\X,$ {where the sets $A^*(x)$ are defined in \eqref{defsetA*}.}  \Halmos\endproof

Analyzing the proofs of Hern\'{a}ndez-Lerma~\cite[Section 4, Theorem]{HL91} and Feinberg and Kasyanov \cite[Theorem~3.3]{MDPSetwise2021}, we obtain the following theorem.

\begin{theorem}\label{thm:ACOI:s} Let Assumptions~S* hold, and let Assumption ${\rm  \underline{B}_{\{\alpha_n\}}}$ hold for a sequence $\{ \a_n\uparrow 1 \}_{n\in\N^*}.$
 Let
\begin{align}
  	u(x) := \liminf_{n\to \infty} u_{\a_n}(x) , \qquad x\in\X .
  	\label{eqn:tu:setwise}
\end{align}
Then there exists a  deterministic  policy $\phi$ satisfying WACOI \eqref{eqn:ACOIOver} with the function $u$ defined in
\eqref{eqn:tu:setwise}.  Therefore, $\phi$ is a deterministic average-cost optimal policy.
\end{theorem}

\proof{Proof.} The proof of optimality inequality \eqref{eq7111} follows the original proof of Hern\'{a}ndez-Lerma~\cite[Section 4, Theorem]{HL91} with minor modifications;  see, also, Feinberg and Kasyanov \cite[Theorem~3.3]{MDPSetwise2021}). Inequality   \eqref{eq7111} implies WACOI \eqref{eqn:ACOIOver} in view of Theorem~\ref{Prop1}.   \Halmos\endproof

The following corollary provides under Assumptions W* or S* sufficient conditions for the validity of ACOI \eqref{eqn:ACOI} under weaker conditions than Assumption~B.
\begin{corollary}\label{corACOIwm}
Let Assumption W* or S* hold, and let there exist a sequence $\{ \a_n\uparrow 1 \}_{n\in\N^*}$ such that Assumption ${\rm  \underline{B}_{\{\alpha_n\}}}$ holds, and $(1-\alpha_n)m_{\alpha_n}\to\underline{w} $ as $n\to\infty.$ Then ACOI \eqref{eqn:ACOI} holds.
\end{corollary}
\proof{Proof.}
Theorems~\ref{thm:ACOI:w} and \ref{thm:ACOI:s} imply that WACOI \eqref{eqn:ACOIOver} holds.  In addition, since $\underline{w}=\overline{w}$, we see that ACOI \eqref{eqn:ACOI} holds.
\Halmos\endproof

Recall the following definitions.

\begin{definition}[Semi-equicontinuity \cite{FKL18a}]\label{def:lowerequicont}
	A sequence $\{f_n\}_{n\in\N^*}$ of real-valued functions on a metric space $(\S,\rho)$
	is called \textit{lower semi-equicontinuous at the point} $s\in \S$ if for each $\ee > 0$
	there exists $\delta>0$ such that
	\begin{equation*}
		 f_n(s^\prime) >  f_n(s) - \ee \qquad 
         \text{ for all } n\in\N^* \text{\ if } \rho(s,s')<{\delta}.
	\end{equation*}
	The sequence $\{f_n\}_{n\in\N^*}$ is called \textit{lower semi-equicontinuous} (\textit{on $\S$}) if it is lower semi-equicontinuous at
	all $s \in \S.$
	A sequence $\{f_n\}_{n\in\N^*}$ of real-valued functions on a metric space $\S$
	is called \textit{upper semi-equicontinuous at the point} $s\in \S$ (\textit{on $\S$})
if the sequence $\{ -f_n \}_{n\in\N^*}$ is lower semi-equicontinuous at the point $s\in \S$ (on~$\S$).
\end{definition}

\begin{definition}[Equicontinuity]\label{def:quicont}
A sequence $\{f_n\}_{n\in\N^*}$  of real-valued functions on a metric space $\S$
is called \textit{equicontinuous at the point} $s\in \S$ (\textit{on $\S$}) if this sequence is
both lower and upper semi-equicontinuous at the point $s\in \S$ (on $\S$).
\end{definition}

The following corollary from Theorem~\ref{thm:ACOI:w} provides a sufficient condition for the validity of ACOI~\eqref{eqn:ACOI} with a relative value function $u$ defined in \eqref{eqn:tu:setwise}.
\begin{corollary}\label{c6.4EF} Let Assumptions~W* hold,  Assumption~${\rm \underline{B}_{\{\alpha_n\}}}$ hold for a sequence $\{ \a_n\uparrow 1 \}_{n\in\N^*}, $ and the sequence of functions $\{u_{\a_n}\}_{n\in\N^*}$ {be} lower semi-equicontinuous.  Then the conclusions of Theorem~\ref{thm:ACOI:w} hold for the function $u$ defined in \eqref{eqn:tu:setwise} for this sequence $\{ \a_n\}_{n\in\N^*}.$
\end{corollary}
\proof{Proof.} Since the sequence of functions $\{u_{\a_n}\}_{n\in\N^*}$ is lower semi-equicontinuous, the functions $u$ defined in \eqref{EQN1} and in \eqref{eqn:tu:setwise} coincide in view of  {\cite[Theorem~3.1(i)]{FKL18a}}.
\hfill\Halmos\endproof

 Consider the following version of the equicontinuity condition (EC) on the discounted relative value functions from \cite{FKL18a}. 

\vspace{0.1in}
\noindent
\textbf{Assumption ${\rm  EC_{\{\alpha_n\}}}$.} { The sequence  of discount factors $\{ \a_n\uparrow 1 \}_{n\in\N^*}$ satisfies the following properties:} \begin{itemize} \item[(i)] the sequence of functions $\{u_{\a_n}\}_{n\in\N^*}$ is equicontinuous; \item[(ii)] there exists a nonnegative measurable function $U(x),$ $x\in\X,$ such that
$U(x)\geq u_{\a_n}(x),$ $n\in\N^*,$ and $\int_{\X} U(y)q(dy|x,a) < +\infty$ for all $x\in\X$ and $a\in \A.$ \end{itemize}

Under each of the Assumptions~W* or \cite[Assumption 4.2.1]{HLL96}, which is stronger than Assumption~S*, and under Assumptions ${\rm  \underline{B}_{\{\alpha_n\}}}$ and ${\rm  EC_{\{\alpha_n\}}}$, there exists a deterministic policy $\phi$ satisfying the average-cost optimality equation (ACOE)
\begin{equation}
    \begin{aligned}
	w^* + u(x) = &c(x,\phi(x)) + \int_{\X} u(y)q(dy|x,\phi(x))\\  =&
  		\min_{a\in \A} \left[ c(x,a) + \int_{\X} u(y)q(dy|x,a) \right],\  x\in\X,
\label{eqn:ACOE}
\end{aligned}
\end{equation}
with $u$ defined in \eqref{EQN1} for the sequence $\{ \a_n\uparrow 1 \}_{n\in\N^*},$ and the function $u$ is continuous; see Feinberg and Liang~\cite[Theorem~3.2]{FLi17} for Assumption~W* and
Hern\'{a}ndez-Lerma and Lasserre~\cite[Theorem 5.5.4]{HLL96}.
We remark that the quantity $w^*$ in \eqref{eqn:ACOE} can be replaced with any other quantity~in~\eqref{eq:7121a}.
In addition, since the {first} equation in \eqref{eqn:ACOE} implies inequality \eqref{eq7111}, every deterministic policy $\phi$ satisfying
\eqref{eqn:ACOE} is average-cost optimal.  Observe that in these cases the function $u$ is continuous (see \cite[Theorem~3.2]{FLi17} for Assumption~W* and \cite[Theorem 5.5.4]{HLL96}), while under conditions of Theorems~\ref{thm:ACOI:w} and \ref{thm:ACOI:s} the corresponding functions $u$ may not be continuous; see Examples~7.1 and 7.2 from \cite{FKL18a}.
Below we provide more general conditions for the validity of the ACOEs.  In particular, under these conditions the relative value functions $u$ may not be continuous.

Now, we introduce Assumption~${\rm  LEC_{\{\alpha_n\}}}$, which is weaker than Assumption~${\rm  EC_{\{\alpha_n\}}}$.  Indeed, Assumption~${\rm  EC_{\{\alpha_n\}}}$(i) is obviously stronger than ${\rm  LEC_{\{\alpha_n\}}}$(i). In view of the Ascoli  theorem (see  \cite[p. 96]{HLL96} or \cite[p. 179]{Roy68}), ${\rm  EC_{\{\alpha_n\}}}$(i) and the first claim in ${\rm  EC_{\{\alpha_n\}}}$(ii) imply ${\rm  LEC_{\{\alpha_n\}}}$(ii). The second claim in ${\rm  EC_{\{\alpha_n\}}}$(ii) implies ${\rm  LEC_{\{\alpha_n\}}}$(iii).
It is shown in Theorem~\ref{thm:acoe:L}  that the ACOEs hold under
Assumptions~W*,  ${\rm  \underline{B}_{\{\alpha_n\}}}$, and ${\rm  LEC_{\{\alpha_n\}}}$.

\vspace{0.1in}
\noindent
\textbf{Assumption ${\rm  LEC_{\{\alpha_n\}}}$.} { The sequence  of discount factors $\{ \a_n\uparrow 1 \}_{n\in\N^*}$ satisfies the following properties:}  \begin{itemize} \item[(i)] the sequence of functions $\{u_{\a_n}\}_{n\in\N^*}$ is lower semi-equicontinuous;
    \item[(ii)] $\lim_{n\to\infty} u_{\a_n} (x)$ exists for each $x\in\X;$ 
\item[(iii)] for each $x\in\X$ and $a\in\A$  the sequence $\{u_{\a_n}\}_{n\in\N^*}$
		is {asymptotically uniformly integrable with respect to the probability measure $q(\,\cdot\,|x,a),$
        that is,
        \[
\lim_{K\to +\infty} \limsup_{n\in \N^*}\int_\X u_{\a_n} (z)q(dz|x,a)=0,\]
which, according to \cite[Theorem 2.2]{FKL18}, is equivalent to the existence of $N\in\N^*$ such that the sequence $\{u_{\alpha_N},u_{\alpha_{N+1}},\ldots \}$ is uniformly integrable with respect to the probability measure $q(\,\cdot\,|x,a).$ }
        \end{itemize}

\begin{theorem}\label{thm:acoe:L}
Let Assumption~W* hold, and let Assumption~ ${\rm  \underline{B}_{\{\alpha_n\}}}$ hold for a sequence $\{ \a_n\uparrow 1 \}_{n\in\N^*} .$  If Assumption ${\rm  LEC_{\{\alpha_n\}}}$ is satisfied for the sequence $\{ \a_n\}_{n\in\N^*},$
then there exists a deterministic  policy $\phi$ such that ACOE \eqref{eqn:ACOE} hold with the function $u(x)$
defined in \eqref{eqn:tu:setwise}.
\end{theorem}

\proof{Proof.}
Since Assumptions~W* and ${\rm  \underline{B}_{\{\alpha_n\}}}$ hold, and $\{ u_{\a_n} \}_{n\in\N^*}$ is lower semi-equicontinuous, then Corollary~\ref{c6.4EF} implies  {the existence of} a deterministic policy $\phi$ satisfying \eqref{eq7111} with $u$ defined in \eqref{eqn:tu:setwise} 
\begin{align}
c(x,\phi(x)) + \int_{\X} u(y)q(dy|x,\phi(x))\le 	w^* + u(x) ,\quad x\in\X .
	\label{eqn:geq}
\end{align}

To prove the ACOE, it remains to prove the opposite inequality to \eqref{eqn:geq}.
According to Feinberg et al.~\cite[Theorem 2(iv)]{FKZ12}, for each $n\in\N^*$ and $x\in\X$ the discounted-cost optimality equation is
$v_{\a_n}(x) = \min_{a\in \A} [ c(x,a) + \a_n \int_{\X} v_{\a_n} (y) q(dy|x,a) ],$  which, by subtracting ${\alpha_n}m_{\a_{n}}$ from both sides and by replacing $\a_n$ with 1, implies that for all $a\in\A$
\begin{align}
	(1-\a_n)m_{\a_n} + u_{\a_n} (x) \leq 
 c(x,a) +  \int_{\X} u_{\a_n} (y) q(dy|x,a),\qquad x\in\X  .
	\label{eqn:transform dcoe}
\end{align}
Let $n\to \infty.$ In view of \eqref{eq:7121a}, Assumptions~${\rm  LEC_{\{\alpha_n\}}}$(ii, iii), and Fatou's lemma~\cite[p. 211]{Shi96},
\eqref{eqn:transform dcoe} imply that for all $a\in\A$
\begin{align}
	w^* + u(x) \leq  c(x,a) + \int_{\X} u(y)q(dy|x,a), \qquad x\in\X .
	\label{eq:acoi:rev}
\end{align}
We remark that the integral in \eqref{eqn:transform dcoe} converges to the integral in \eqref{eq:acoi:rev} since the sequence $\{u_{\a_n}\}_{n\in\N^*}$ converges pointwise to $u$ and is u.i.; see  { \cite[Theorem~2.1]{FKL18a}}.
Then, \eqref{eq:acoi:rev} implies
\begin{align}
	w^* + u(x) \leq \min_{a\in \A} [ c(x,a) + \int_{\X}                                                                                                       u(y)q(dy|x,a)]\leq c(x,\phi(x)) + \int_{\X} u (y) q(dy|x,\phi(x)),\quad x\in\X .
	\label{eqn:leq}
\end{align}
 Thus, \eqref{eqn:geq} and \eqref{eqn:leq} imply
\eqref{eqn:ACOE}.
\hfill\Halmos\endproof

In the following example, Assumptions~W*, ${\rm  \underline{B}_{\{\alpha_n\}}}$, and ${\rm  LEC_{\{\alpha_n\}}}$ hold. Hence the ACOEs hold. However, Assumption~${\rm EC_{\{\alpha_n\}}}$ does not hold.
Therefore, Assumption~${\rm  LEC_{\{\alpha_n\}}}$ is more general than Assumption~${\rm  EC_{\{\alpha_n\}}}$.

\begin{example}{\rm (\cite[Example~7.1]{FKL18a})}\label{ex:mdpLEC}
{\rm
	Consider $\X = [0,1]$ equipped with the Euclidean metric
	and $\A = \{ a^{(1)} \} .$ The transition probabilities are
	$q(0|x,a^{(1)}) = 1$ for all $x\in\X.$ The cost function is $c(x,a^{(1)}) = \h \{ x \neq 0 \},$ $x\in\X.$ Then the discounted-cost value is $v_\a (x) = u_\a (x) =\h \{ x \neq 0 \},$  $\a\in[0,1)$ and $x\in\X,$ and the average-cost value is $w^* = w (x) = 0,$ $x\in\X.$
	It is straightforward to see that Assumptions~W* and ${\rm  \underline{B}_{\{\alpha_n\}}}$ hold.
	In addition, since the function $u(x) = \h \{ x \neq 0 \}$ is lower semi-continuous, but it is not continuous,
	the sequence of functions $\{u_{\a_n}\}_{n\in\N^*}$ is lower semi-equicontinuous,  but it is not
	equicontinuous for each sequence $\{ \a_n\uparrow1 \}_{n\in\N^*}.$
	Therefore, Assumption~${\rm  LEC_{\{\alpha_n\}}}$ holds since $0\leq u_{\a_n} (x)\leq 1,$ $x\in\X,$
	and Assumption~${\rm  EC_{\{\alpha_n\}}}$ does not hold.
	The \eqref{eqn:ACOE} holds with $w^* = 0,$ $u (x) = \h \{ x \neq 0 \},$ and
	$\phi (x) = a^{(1)},$ $x\in\X.$
	\hfill\Halmos\endproof
}
\end{example}

The following theorem states the validity of the ACOE under
Assumptions~S*, ${\rm  \underline{B}_{\{\alpha_n\}}}$, and ${\rm  LEC_{\{\alpha_n\}}}${\rm(ii,iii)}.

\begin{theorem}\label{thm:acoe:S}
Let Assumption~{S}* hold, and let Assumption~ ${\rm  \underline{B}_{\{\alpha_n\}}}$ hold for a sequence $\{ \a_n\uparrow 1 \}_{n\in\N^*} .$ If Assumptions~${\rm  LEC_{\{\alpha_n\}}}${\rm(ii,iii)} are satisfied for the sequence $\{ \a_n\}_{n\in\N^*},$
then there exists a deterministic  policy $\phi$ such that ACOE \eqref{eqn:ACOE} holds with the function $u(x)$
defined in \eqref{eqn:tu:setwise}.
\end{theorem}

\proof{Proof.}
According to Theorem~\ref{thm:ACOI:s},
if Assumptions~S* and ${\rm  \underline{B}_{\{\alpha_n\}}}$ hold, then we have that: (i) equalities in \eqref{eq:7121a} hold;
(ii) there exists a deterministic policy $\phi$ satisfying  ACOI \eqref{eqn:geq} with the function $ u$ defined
in \eqref{eqn:tu:setwise}; and
(iii) for each $n\in\N^*$ and $x\in\X$ the discounted-cost optimality equation is
$v_{\a_n}(x) = \min_{a\in \A} [ c(x,a) + \a_n \int_{\X} v_{\a_n} (y) q(dy|x,a) ].$
Therefore, the same arguments as in the proof of Theorem~\ref{thm:acoe:L} starting from \eqref{eqn:transform dcoe}
imply the validity of \eqref{eqn:ACOE} with $u$ defined
in \eqref{eqn:tu:setwise}.
\hfill\Halmos\endproof

Observe that the MDP described in Example~\ref{ex:mdpLEC}
also satisfies Assumptions~S*, ${\rm \underline{B}_{\{\alpha_n\}}}$, and ${\rm  LEC_{\{\alpha_n\}}}${\rm(ii,iii)}.
We provide Example~\ref{ex:mdpLEC:s},
in which Assumptions~S*, ${\rm  \underline{B}_{\{\alpha_n\}}}$, and ${\rm  LEC_{\{\alpha_n\}}}${\rm(ii,iii)} hold.
Hence, the ACOEs hold. However, Assumptions~W*, ${\rm  LEC_{\{\alpha_n\}}}$(i), and ${\rm  EC_{\{\alpha_n\}}}$ do not hold.

\begin{example}{\rm (\cite[Example~7.2]{FKL18a})}\label{ex:mdpLEC:s}
{\rm
	Let $\X = [0,1]$ and $\A = \{ a^{(1)} \} .$ The transition probabilities are
	$q(0|x,a^{(1)}) = 1$ for all $x\in\X.$ The cost function is $c(x,a^{(1)}) = D(x),$ where $D$ is the
	Dirichlet function defined as
	\begin{align*}
		D (x) =
		\begin{cases}
			0, & \text{if } x \text{ is rational,} \\
			1, & \text{if } x \text{ is irrational,}
		\end{cases}
		\qquad x\in\X.
	\end{align*}
	Since there is only one available action, Assumption~S* holds.
	The discounted-cost value is $v_\a (x) = u_\a (x) = D(x)=u(x),$  $\a\in[0,1)$ and $x\in\X,$
	and the average-cost value is $w^* = w (x) = 0,$ $x\in\X.$ Then Assumptions~${\rm \bf \underline{B}_{\{\alpha_n\}}}$ and ${\rm  LEC_{\{\alpha_n\}}}${\rm(ii,iii)} hold.
	Hence, the ACOEs \eqref{eqn:ACOE} hold with $w^* = 0,$ $                                                                                                                      u (x) = D(x),$ and $\phi (x) = a^{(1)},$ $x\in\X.$  Thus, the average-cost relative function $u$ is not lower semi-continuous.
	However, since the function $c(x,a^{(1)})=D(x)$ is not lower semi-continuous, Assumption~W* does not hold.
	 Since the function $u(x) =u_\a (x) = D(x)$ is not lower semi-continuous,  Assumptions~${\rm  LEC_{\{\alpha_n\}}}$(i) and ${\rm  EC_{\{\alpha_n\}}}$ do not hold either.
	\hfill\Halmos\endproof
}
\end{example}

We recall that, in view of Theorem~\ref{Prop1}, $w^*=\overline{w}$ under assumptions of {this theorem.} The following theorem provides sufficient conditions for $w^*=\overline{w}=\underline{w}.$  While ACOE \eqref{eqn:ACOE} is a stronger fact than WACOE \eqref{eqn:ACOIOver}, Corollary~\ref{ceqACEstr} provides sufficient conditions for the optimality equality which is stronger than ACOE \ref{eqn:ACOI}.

\begin{corollary} \label{ceqACEstr} Let  assumptions of either Theorem~\ref{thm:acoe:L} or Theorem~\ref{thm:acoe:S} hold. If, in addition
$\lim_{n\to\infty}(1-\alpha_n)m_{\alpha_n}=\underline{w},$ then $w^*=\overline{w}=\underline{w}$ and ACOE \eqref{eqn:ACOE} holds with $w^*$ substituted with $\underline{w}.$
\end{corollary}
\proof{Proof.} The proof follows from the arguments provided after the formulation of Theorem~\ref{Prop1}.~\hfill\Halmos\endproof

We remark that \cite[Example 4.1]{MDPSetwise2021} satisfies the assumptions of Corollary~\ref{ceqACEstr}, but it does not satisfy Assumption B.

\vskip 6mm
\noindent{\bf Acknowledgments}

\noindent
This research was partially supported by SUNY System Administration under SUNY Research Seed Grant
Award 231087 and by the U.S. Office of Naval Research (ONR) under Grants N000142412608 and N000142412646.

\end{document}